\documentclass{article}
\usepackage[margin=1in]{geometry}

\usepackage[backend=biber,style=alphabetic,sorting=nyt]{biblatex}
\usepackage{amssymb,amsfonts,amsmath,bm}

\usepackage{booktabs,array}
\newcolumntype{L}[1]{>{\raggedright\arraybackslash}p{#1}}

\usepackage{xcolor}
\definecolor{darkgreen}{rgb}{0.0, 0.5, 0.0}

\usepackage{graphicx}

\usepackage[
    colorlinks=true,
    linkcolor=purple,
    citecolor=darkgreen,
    urlcolor=red
]{hyperref}

\usepackage{caption}
\usepackage{todonotes}
\usepackage{lineno}

\title{Scientific Applications Leveraging Randomized Linear Algebra}
\author{
Vivak Patel\thanks{Statistics at the University of Wisconsin, Madison. Supported by NSF 2309445.}\and 
D. Adrian Maldonado\thanks{Mathematics and Computer Science at Argonne National Laboratory. Supported by NSF 2309446.}\and
Maksim Melnichenko\thanks{Innovative Computing Laboratory at the University of Tennessee, Knoxville. Supported by NSF 2004541.}\and 
Nathaniel Pritchard\thanks{Mathematical Institute at the University of Oxford.}\and 
Vishwas Rao\thanks{Mathematics and Computer Science at Argonne National Laboratory.}\and 
Elizaveta Rebrova\thanks{Operations Research and Financial Engineering at Princeton University.}\and 
Sriram Sankararaman\thanks{Computer Science, Human Genetics, Computational Medicine at the University of California, Los Angeles. Supported by NIH R35GM153406 and NSF 1943497.}\and
Marcel Schweitzer\thanks{School of Mathematics and Natural Sciences, Bergische Universität Wuppertal, Wuppertal, Germany.}
}
\date{\today}

\newcommand{\vct}[1]{\bm{#1}} 
\newcommand{\mat}[1]{\bm{#1}} 
\newcommand{\R}{\mathbb{R}} 
\newcommand{\C}{\mathbb{C}} 
\newcommand{\bigo}[1]{\mathcal{O}(#1)} 

\newcommand{\tr}{\mathrm{tr}}
\newcommand{\herm}{^{\dagger}}

\begin{document}

\maketitle

\begin{abstract}
This report showcases the role of, and future directions for, the field of Randomized Numerical Linear Algebra (RNLA) in a selection of scientific applications.
These applications span the domains of imaging, genomics and dynamical systems, and are thematically connected by needing to perform linear algebra routines on large-scale matrices (with up to quantillions of entries).
At such scales, the linear algebra routines face typical bottlenecks: memory constraints, data access latencies, and substantial floating-point operation costs.
RNLA routines are discussed at a high-level to demonstrate how these routines are able to solve the challenges faced by traditional linear algebra routines, and, consequently, address the computational problem posed in the underlying application.
For each application, RNLA's open challenges and possible future directions are also presented, which broadly fall into the categories: creating structure-aware RNLA algorithms; co-designing RNLA algorithms with hardware and mixed-precision considerations; and advancing modular, composable software infrastructure. 
Ultimately, this report serves two purposes: it invites domain scientists to engage with RNLA; and it offers a guide for future RNLA research grounded in real applications.
\end{abstract}

\tableofcontents

\newpage

\pagebreak

\section{Overview}

Randomized numerical linear algebra (RNLA) uses randomization to enhance classical numerical linear algebra methods by reducing their computational costs, accelerating their performance, and improving their numerical stability.
For over a decade, RNLA has enjoyed strong theoretical foundation and numerical experiments that consistently demonstrates its effectiveness \cite[e.g.,][]{halko2011finding,mahoney2011random}.%
\footnote{For readers interested in learning about RNLA, an incomplete set of starting points include \cite{halko2011finding,woodruff2014sketching,drineas_randnla_2016,drineas2017lectures,martinsson2020randomized,murray_randomized_2023}.}
Of equal importance, RNLA is demonstrating its impact on real applications. 

A selection of applications are showcased in this paper with two specific goals. 
The first is to expose domain scientists to the utility of RNLA with real use cases, which we hope serves as an invitation for these scholars to engage with RNLA.
The second is to collate open challenges in RNLA theory, methodology and software to guide future developments in the field.

Although these applications are from a sampling of RNLA researchers (biased towards those attending the 2024 SIAM Conference on Applied Linear Algebra in Paris, France),
the applications, the underlying RNLA methodologies, and their future directions have coalesced around a few directions:
the applications are data intensive;
the data or resulting matrices of interest have informative compressions constructed via random matrix multiplication;
and application-motivated future directions center around useful software that allows for adapting RNLA methods to specific problem structure and hardware resources.

\subsubsection*{Data-intensive Applications}

The application domains presented below fall into three areas: 
imaging (see \S\ref{application:ct} and \S\ref{application:hyperspectralmixing}); 
genomics (see \S\ref{application:gwas1}, \S\ref{application:gwas2} and \S\ref{application:snp}); 
and dynamical systems (see \S\ref{application:oi} to \S\ref{application:qcd_trace}).
Though seemingly disparate, these application domains share a common feature: the computational problems emerging in these domains are generally data-intensive, requiring the application of numerical linear algebra methods to matrices at large scales (see Table \ref{table:matrix_scales}).
Owing to these scales, numerical linear algebra methods are often constrained by the memory available in hardware, the cost of moving information, or the cost of performing arithmetic operations on the underlying matrices. 
In the applications discussed below, these computational bottlenecks are addressed by using randomized techniques to compress---that is, to create smaller representations that preserve the essential information needed for computation---either the underlying data or the resulting matrices, and
then applying numerical linear algebra algorithms to these compressed matrices. 

\begin{table}[hb!]
\centering
\caption{Matrix dimensions for solving example computational problems in the listed application area.
\label{table:matrix_scales}}
\begin{tabular}{L{1.75in}L{1.75in}} \toprule
    \textbf{Application Area} & \textbf{Matrix Dimensions} \\ \midrule 
    Imaging & $\bigo{10^6 \times 10^6}$ \\
    Genomics & $\bigo{10^6 \times 10^8}$ \\
    Dynamical Systems & $\bigo{10^8 \times 10^9}$ \\ \bottomrule 
\end{tabular}
\end{table}

\subsubsection*{The Role of Randomized Numerical Linear Algebra}

One of the major insights of randomized numerical linear algebra exploited in the applications below is: large matrices frequently encountered in practice have highly informative compressions (e.g., low-rank projections), which can be inexpensively approximated using multiplication with a random matrix.%
\footnote{As seen below, there is usually some additional numerical linear algebra calculation performed, but it is relatively inexpensive to perform on the compression.}
The compressions can then be used to find an approximate solution of a sub-problem in the target application's computational workflow, or can be iterated to deliver a more accurate solution. 

The compressions are used for distinct linear algebra tasks in the applications listed below: for solving large-scale least squares problems; for generating informative matrix decompositions; or for approximating matrix functions
(see Table \ref{table:tasks-applications}).
Importantly, the linear algebra tasks are not fixed to an application domain; for example, least squares problems appear in applications arising in imaging, genomics and dynamical systems. 
That is, the range of RNLA techniques seems to be useful across a number of scientific domains.
Thus, the main focus of RNLA in the applications below is to ensure that these
methods are adapted to the needs of the specific computational problem.

\begin{table}[ht!]
\centering
\caption{Example numerical linear algebra tasks and the applications that correspond to these tasks.}
\label{table:tasks-applications}
\begin{tabular}{L{2in}L{2in}} \toprule 
    \textbf{Task} & \textbf{Applications} \\ \midrule 
    Least Squares Problems & \S\ref{application:ct}, \S\ref{application:gwas1}, \S\ref{application:oi} and \S\ref{application:da} \\
    Matrix Decompositions & \S\ref{application:hyperspectralmixing} and \S\ref{application:snp} \\ 
    Matrix Functions & \S\ref{application:gwas2}, \S\ref{application:qcd_overlap} and \S\ref{application:qcd_trace} \\\bottomrule 
\end{tabular}
\end{table}

\subsubsection*{Application-motivated Future Directions}

While advancements in RNLA continue to be driven by researchers in the field,
RNLA's advancements are also being driven by the needs of applications.
Many of these needs resonate in the future directions and open challenges discussed in each application below.
Some themes are summarized here and are elaborated in \S\ref{section:conclusion}.

\begin{itemize}
    \item {\it Structure-aware algorithms.} Future work must focus on developing adaptive, randomized methods that automatically detect and preserve problem-specific structures, ensuring domain-specific accuracy.
    \item {\it Hardware and precision co-design.} RNLA must be optimized for modern hardware through cache-aware and GPU-efficient algorithms, and enhanced with precision-aware strategies that balance performance with theoretical rigor.
    \item {\it Software infrastructure.} The transition from proof-of-concept RNLA implementations to production-ready software libraries remains a critical bottleneck for widespread adoption. This includes developing accessible interfaces, mature abstractions, reusable data structures, and modular software that integrate seamlessly with application-specific workflows and data formats.
\end{itemize}

\noindent We now turn our attention to the specific applications.

\newpage 
\section{Applications}
\subsection{Computed Tomography} \label{application:ct}

\begin{figure}[h!]
    \centering
    \includegraphics[width=0.5\linewidth]{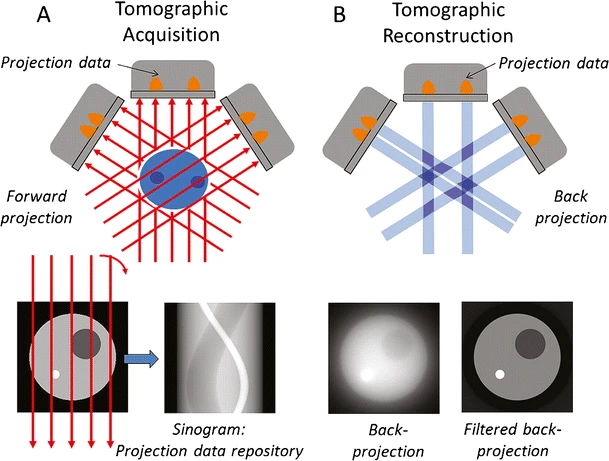}
    \caption{A demonstration of the CT process from imaging (A) to reconstruction via two example algorithms (B). The imaging generates a sinogram, which is then algorithmically processed to approximate the original image. This is a reproduction of  Figure 1 from \cite{seibert2014iterative} with permission from the authors.}
    \label{fig:ct}
\end{figure}

Computed Tomography (CT) is a non-invasive imaging technique that uses X-rays to obtain images of the body. This technique reconstructs an image from several X-ray measurements taken at multiple angles and, as a consequence, it is able to detect medical issues that would not be captured from standard X-ray scans (see Figure \ref{fig:ct}).  This image reconstruction is one of the fundamental challenges of CT. Mathematically, it can be formulated as the solution of a large-scale overdetermined system \cite{CTBook}.

\subsubsection*{High-Level Mathematical Description}
CT reconstruction begins by encapsulating the interactions of X-ray beams with an object of interest through a matrix $\mat{A} \in \R^{m \times n}$, where entry $\mat{A}_{ij}$ encodes the length of the $i^\mathrm{th}$ X-ray beam's interaction with the $j^\mathrm{th}$ position of the object; 
and a vector $\vct{b} \in \R^{m}$ contains the modulated X-ray intensities picked up by the sensor after passing through the object. 
CT reconstruction solves for the vector $\vct{x} \in \R^n$, which are the attenuation coefficients of the object corresponding to its material properties, from the linear system $\mat{A}\vct{x} = \vct{b}$.

Since modern CT systems acquire many more measurements than the number of unknowns (\(m > n\)), reconstruction is done by solving,
\begin{equation}
\min_{\vct{x} \in \R^n} \| \mat{A}\vct{x} - \vct{b} \|_2.
\end{equation}
The fundamental challenge arises because \(\mat{A}\) can become extremely large as the resolution grows and number of beams increase. For instance, reconstructing a \(1024 \times 1024\) pixel image leads to \(n = 10^6\) unknowns, and \(m\) can be on the order of several million, making standard direct methods for solving large systems (e.g., factorizing \(\mat{A}\)) impractical in terms of both storage and computation.

\subsubsection*{The Potential Role of RNLA}
The Kaczmarz method, known as the Algebraic Reconstruction Technique (ART) in imaging applications, has been a foundational approach in computed tomography reconstruction. Its iterative row projection strategy made it one of the earliest successful methods for CT reconstruction, particularly appealing because of its simple implementation and ability to handle the large, sparse systems typical in medical imaging. 
However, its deterministic nature limits its effectiveness in handling the inherent redundancies in CT data, which arise 
from similarly oriented X-ray beams (or small shifts in detector positions) yielding highly correlated rows in $\mat{A}$.

The randomized Kaczmarz method, introduced by Strohmer and Vershynin, represented a significant advancement to address this challenge \cite{Strohmer2008}. By randomly sampling from the available measurements, the randomized Kaczmarz method can better avoid the redundancies that limited the cyclic approach of ART. Furthermore, the randomized Kaczmarz method enjoys an exponential convergence rate (in expectation)---in stark contrast to the convergence rate of ART, which depends on the information dissimilarity between sequential equations.

\subsubsection*{Outlook and Open Challenges}
The next frontier lies in developing more sophisticated randomized techniques that exploit the specific mathematical structure of CT problems. Block methods that process multiple projections simultaneously using randomized dimension reduction could leverage the inherent redundancy in CT measurements \cite{Needell2015, Patel2023block}; adaptive sampling strategies could potentially accelerate convergence by focusing computational effort on the most informative measurements \cite{Gower2021, Patel2023block}; and robustness strategies can reduce the impact of potential corruptions or errors in measurements \cite{haddock2022quantile, cheng2022block}.

\subsection{Hyperspectral Unmixing} \label{application:hyperspectralmixing}



\begin{figure}[h!]
    \centering
    \includegraphics[width=0.5\linewidth]{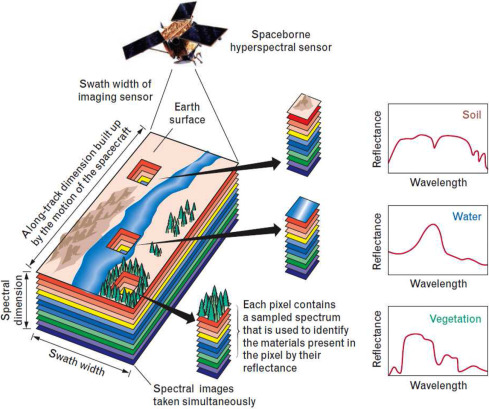}
    \caption{A graphic representation of hyperspectral unmixing. This is a reproduction of Figure 1 from \cite{dobigeon2016linear} with permission from the authors. This image originally appeared in \cite{shaw2003spectral}.
    \label{fig:hm}}
\end{figure}

Hyperspectral images are generated by specialized cameras that record the intensity of light reflected by real objects at many different wavelengths (usually, several hundred). 
Hyperspectral images are often used to identify the different materials ({endmembers}) composing the object being captured and their proportions (abundances), such as in identifying mineral compositions of a hyperspectral image of a geographic area \cite{clark1992mapping}. The process of recovering endmembers and abundances from hyperspectral images is called (blind) hyperspectral unmixing (see Figure \ref{fig:hm}), and can be mathematically encapsulated as a matrix decomposition \cite{hyperspectral, gillis2020nonnegative}.

\subsubsection*{High-level Mathematical Description}

Consider a hyperspectral image of $m$ pixels, where each pixel 
corresponds to non-negative intensity measurements of $n$ different wavelengths. 
A hyperspectral image can be represented by a matrix $\mat{A} \in \R^{m \times n}$, where $\mat{A}_{ij}$ is the intensity measurement of wavelength $j$ at pixel $i$. 

In some cases, the hyperspectral image can be assumed to be a linear combination of spectral signatures of its endmembers.
In this case, the hyperspectral image's endmembers and their abundances at each pixel can be recovered by solving 
\begin{equation}\label{eq:nmf}
    \min_{\mat{U},\mat{V}\geq 0} \|\mat{A}-\mat{U}\mat{V}^T\|_F^2,
\end{equation}
where $\mat{U} \in \R^{m \times r}$ has non-negative entries representing the abundances of $r$ endmembers at each pixel; and $\mat{V} \in \R^{n \times r}$ has non-negative entries representing the spectral signature of the $r$ endmembers. 

The problem in \eqref{eq:nmf}, known as non-negative matrix factorization (NMF), is (in general) NP-hard \cite{nmf_np_hard} to solve and might have multiple solutions.
Despite this, NMF problems can often be solved by iterative algorithms, including multiplicative updates and (hierarchical) alternating least squares.
However, when the dimensions of $\mat{A}$ are large, even storing the matrix $\mat{A}$ becomes prohibitive, leaving alone the infeasible memory movement costs of applying iterative algorithms for solving the NMF problem. 

\subsubsection*{The Potential Role of RNLA}

This memory challenge motivates the following approach to the NMF problem: compute a linear function of $\mat{A}$, 
$\mathcal{L}(\mat{A})$, that reduces the memory burden of working with $\mat{A}$, and then use this compression to (approximately) solve the NMF problem. 
One successful choice for this approach is to use a random projection matrix followed by some approximate decomposition technique.
Examples in this category include random projection streams (iterative random sketching approach that requires multiple sketches and passes over $\mat{A}$) \cite{yahaya2021random}; data-adapted sketching with randomized SVD techniques \cite{halko2011finding,tepper2016compressed,erichson2018randomized,hayashi2024randomized}; and alternative formulations of the compressed NMF problem (including efficient regularizers) that have provably comparable optimal solutions to the original problem's solutions \cite{chaudhry2024learning}.

\subsubsection*{Outlook and Open Challenges}
For two-dimensional images, one challenge is to developed specific RNLA methods that can be used in a hierarchical or sparse context \cite{gillis2020nonnegative}.
For three-dimensional or temporal scenes that give rise to information stored as a tensor, the primary challenge is to develop appropriate extensions of the aforementioned two-dimensional techniques to these higher-dimensional cases.

\subsection{Genome-wide Association Studies I: Linear Models} \label{application:gwas1}

\begin{figure}[ht!]
    \centering
    \includegraphics[width=0.75\linewidth]{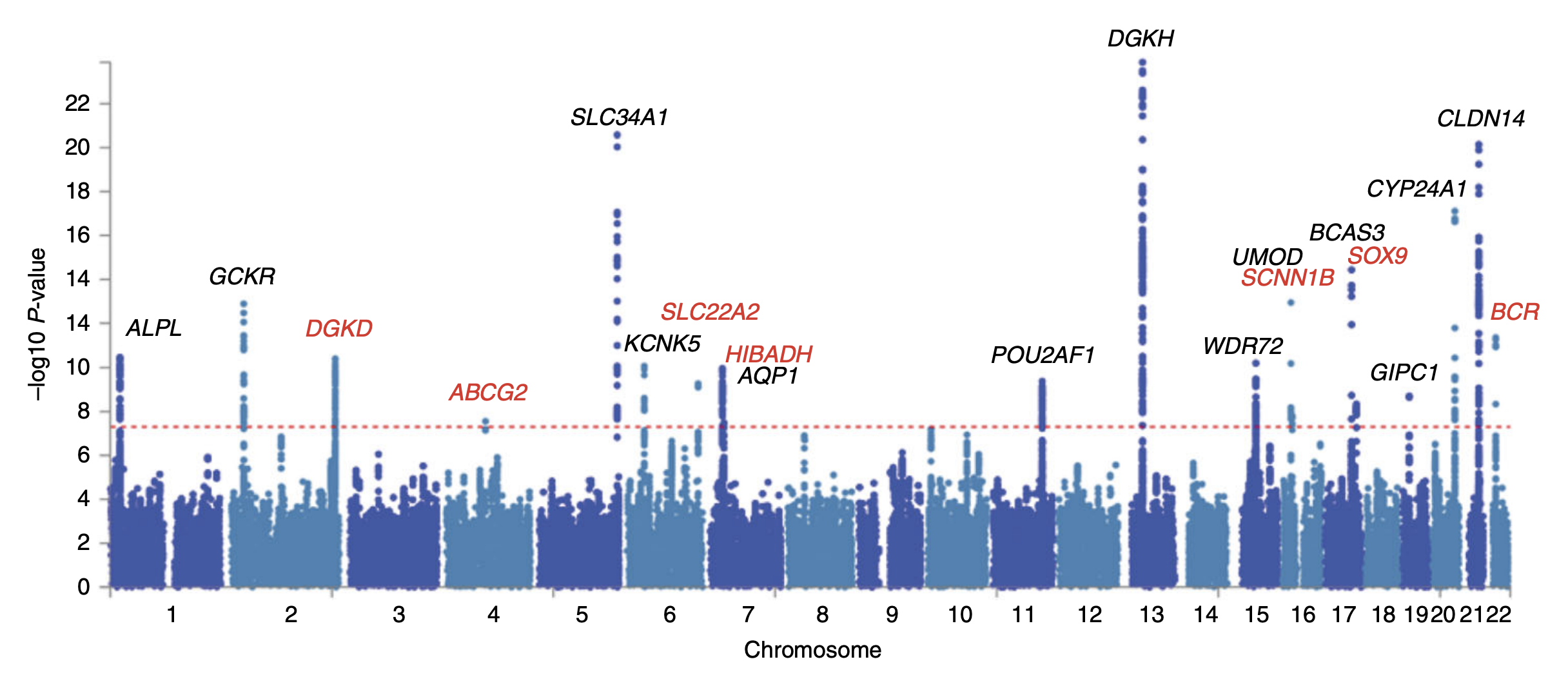}
    \caption{An example in which different areas of the genome (on the x-axis) are associated with a specific kidney disease. The association level is indicated by the y-axis: the higher the value, the more likely the SNP is associated with the disease. This is Figure 1b from \cite{howles2019genetic} and is reproduced under CC-BY 4.0. 
    \label{fig:gwas}}
\end{figure}

Genome-wide Association Studies (GWASs) are tools used to explore the relationship between phenotypes (e.g., a disease status) and genetic information \parencite{uffelmann2021Genomewidea}. 
In particular, GWASs attempt to identify how important certain genetic markers, called single nucleotide polymorphisms (SNPs), are in discerning between phenotypes (see Figure \ref{fig:gwas}). As one option, a GWAS can be formulated mathematically as solving a regularized least squares problem.  

\subsubsection*{High-Level Mathematical Description}
GWASs typically proceed by collecting a set of $m$ individuals and, for each person $i \in \lbrace 1,\ldots,m \rbrace$, measure some phenotype about them, $b_i \in \R$, 
and genetic information about them, $\vct{a}_i \in \R^n$.
A GWAS assumes that there is some unknown linear combination of the genetic information, $\vct{a}_i^\intercal \vct{x}$, that is informative about the phenotype, $b_i$; and attempts to compute the weights of the linear combination, $\vct{x} \in \R^n$. 

A standard approach is to use linear regression to estimate $\vct{x}$, yet this results in poor predictive performance as the resulting system is under-determined: $n$ (the number of genetic predictors) is typically on the order of $10^5$ to $10^8$, while $m$ (the number of individuals) is on the order of $10^3$ to $10^6$ \parencite{devlaming2015Current,ukb.nature.2018,mvp,aou}.
To adjust for this, linear regression can employ a penalty---such as an $\ell_1$ (LASSO) penalty or an $\ell_2$ penalty (Ridge Regression) on $\vct x$---that emphasizes a parsimonious model.
While the selection properties of an $\ell_1$ penalty are appealing, linear regression with an $\ell_2$ penalty can have superior predictive performance for certain GWASs
\parencite{frank1993Statistical,bovelstad2007Predicting,vanwieringen2009Survival}.

Linear regression with an $\ell_2$ penalty, better known as Ridge Regression, requires solving 
\begin{equation} \label{eqn-gwas-ridgeregression}
    \min_{\vct{x} \in \R^n} \Vert \mat{A}\vct{x} - \vct{b} \Vert_2^2 + \lambda \Vert \vct{x} \Vert_2^2,
\end{equation}
where the $i^\mathrm{th}$ row of $\mat{A}\in \R^{m \times n}$ is $\vct{a}_i$; the $i^\mathrm{th}$ entry of $\vct{b} \in \R^m$ is $b_i$; and  $\lambda > 0$ is user-selected penalty parameter \parencite{whittaker2000Markerassisted}.%
\footnote{Choosing the penalty parameter is based on cross-validation methods, or intimate knowledge of the data and problem. Furthermore, Ridge Regression for GWAS can have certain variations where a different penalty parameter is applied to each component of $\vct{x}$.}
The Ridge Regression problem's solution is explicitly given by 
\begin{equation} \label{eqn-gwas-ridgeregression-explicit}
    \vct{x} = (\mat{A}^\intercal \mat{A} + \lambda \mat{I})^{-1}\mat{A}^\intercal \vct{b},
\end{equation}
where $\mat{I}$ is the $n \times n$ identity matrix.

\subsubsection*{The Potential Role of RNLA}
Owing to the size of $\mat{A}$ (on the order of $10^3 \times 10^5$ to $10^6 \times 10^8$), solving (\ref{eqn-gwas-ridgeregression}) or directly computing (\ref{eqn-gwas-ridgeregression-explicit}) is nontrivial. 
While certain iterative methods can be used for this purpose, RNLA offers a potentially better solution.

To be specific, let $\vct{z} \in \R^n$ be a random variable with expectation zero and covariance matrix $\lambda \mat{I} \in \R^{n \times n}$. Let $\mat{Z} \in \R^{m \times n}$ be a matrix whose rows are independent copies of $\vct{z}$.
Then, the Ridge Regression problem is equivalent to solving
\begin{equation} \label{eqn-gaws-ridgeregression-reform}
    \min_{\vct{x} \in \R^n} \mathbb{E} \left[ \Vert (\mat{A} + \mat{Z})\vct{x} - \vct{b} \Vert_2^2 \right].
\end{equation}
That is, this formulation of the Ridge Regression problem can be successfully addressed with a variety of randomized (streaming) linear solvers that are designed for such large-scale problems \parencite{patel2023Randomized}.

\subsubsection*{Outlook and Open Challenges}
While the alternative formulation of the Ridge Regression problem can be solved using randomized linear solvers, whose theory is developed in \parencite{patel2023Randomized},
there are a number of methodological and software questions.
Methodologically, how can we design $\mat{Z}$ to respect the sparsity in $\mat{A}$ in order to reduce computational costs? Furthermore, how can we use the solution of one Ridge Regression problem with one choice of penalty parameter to inform the solution of the same Ridge Regression problem with a different penalty parameter? 
From a software perspective, how can we work with geneticists to create a software that implements RNLA methodologies using the mature data formats and standards that are common in the genetics discipline?
By addressing these different questions, we can help geneticists unlock more efficient tools to support the critical exploratory work done in GWASs.

\subsection{Genome-wide Association Studies II: Linear Mixed-Effect Models} \label{application:gwas2}

In a GWAS, the most common approach for identifying important genetic regions that discern phenotypes is to model the relationship between the genes and phenotypes using linear mixed-effect models \cite{yang.ng.2010,yang.ajhg.2011}.
While many methods exist for estimating the resulting model,
one approach is based on the method of moments (MoMs): a method that requires computing the traces of first and second powers of large matrices.

\subsubsection*{High-Level Mathematical Description}
Let $\vct{b} \in \R^m$ represent the phenotypes of $m$ individuals, and 
let $\mat{A} \in \R^{m \times n}$ encapsulate $n$ pieces of genetic information (i.e., SNPs) for the $m$ individuals. 
In a linear mixed-effect model, 
$\vct{b}$ is assumed to be realized from a normal distribution with mean zero and variance 
\begin{equation}
    \mat{K}(\sigma_g^2, \sigma_e^2) = \frac{\sigma_g^2}{m} \mat{A} \mat{A}^\intercal + \sigma_e^2 \mat{I},
\end{equation}
where $\mat{I}$ is the $m \times m$ identity matrix; and $\sigma_g^2, \sigma_e^2 > 0$ are unknown parameters that are used to interpret how well a phenotype can be explained by genetic information.

  The values of $\sigma_g^2$ and  $\sigma_e^2$ can be estimated using classical methods such as (restricted) likelihood techniques \parencite{reml}. 
However, this approach to estimation requires optimization methods that must repeatedly calculate $\mat{K}$, and compute its determinant and inverse (as applied to vectors), which become infeasible for the relevant scales of $m$ (i.e., $10^3$ to $10^6$) and $n$ (i.e., $10^5$ to $10^8$). 
An alternative estimation approach is to use {method of moments} (MoM) \parencite{MoM}, which requires repeatedly solving
\begin{gather}
    \begin{bmatrix}
        tr({\mat{K}}^2) & tr({\mat{K}})\\
        tr({\mat{K}}) & m
    \end{bmatrix}
        \begin{bmatrix}
        \sigma^2_g \\
        \sigma^2_e
        \end{bmatrix}
    =
            \begin{bmatrix}
                {\vct{b}}^{\intercal}{\mat{K}}{\vct{b}}\\
                {\vct{b}}^{\intercal}{\vct{b}},
            \end{bmatrix} \label{eqn:mom}
\end{gather}
where $\mat{K}$ is the value $\mat{K}(\sigma_g^2, \sigma_e^2)$ for iteratively updated values of $\sigma_g^2$ and $\sigma_e^2$.

\subsubsection*{The Potential Role of RNLA}
While the system in \eqref{eqn:mom} is $2 \times 2$ and is trivial to solve, computing the coefficients of this system requires calculating $\mat{K}, \mat{K}^2 \in \R^{m \times m}$ and computing their traces. For the relevant scales of $m$, this can be prohibitively expensive.

Fortunately, RNLA methods, such as unbiased randomized trace estimators, can be used that only require matrix-vector products between $\mat{A}$ and $s$ zero-mean, identity-covariance random vectors \parencite{Avron11,Wu19,pazokitoroudi.ncomm.2020}. 
The resulting randomized algorithm, the randomized method of moments (RMoM), has a substantially discounted complexity of $\bigo{nms}$.
Empirical results suggest that the RMoM achieves statistical and computational efficiency using small values of $s$ (e.g., $10$ to $100$) \parencite{Wu19,pazokitoroudi.ncomm.2020}.

\subsubsection*{Outlook and Open Challenges}
While prior works \cite{Wu19, pazokitoroudi.ncomm.2020} have empirically shown that the RMoM estimator is computationally efficient while losing relatively little statistical efficiency, it is important to investigate the characteristics of the data that influence statistical and computational efficiency and how choices of the randomized trace approximation can be leveraged to yield estimators that are both computationally and statistically efficient. From a practical perspective, it would be useful to have adaptive approaches to choose the number of random vectors $s$ to be able to achieve specific levels of statistical and computational efficiency.

\subsection{Analyzing Single Nucleotide Polymorphisms} \label{application:snp}

\begin{figure}[!ht]
    \centering
    \includegraphics[width=0.7\linewidth, trim=1cm 15cm 0cm 0cm, clip]{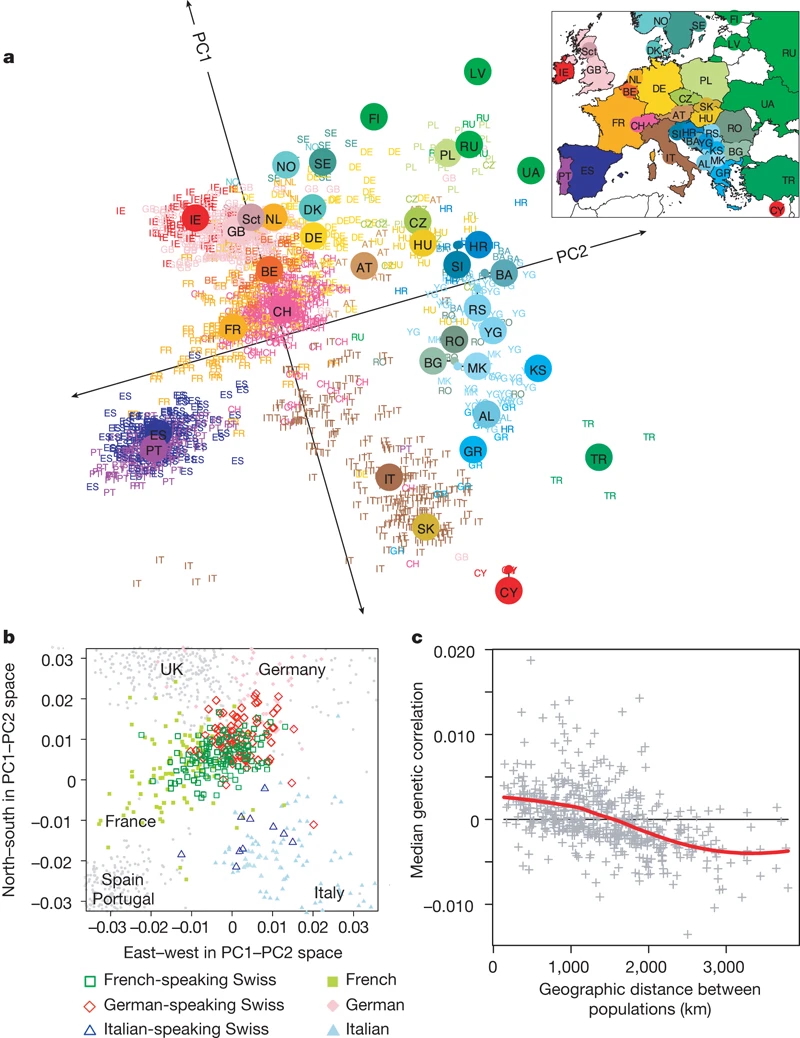}
    \caption{An example analysis of SNP data collected from 30,000 Europeans, which demonstrates how the genetic information mimics geographic regions of Europe. This is Figure 1a from \cite{novembre2008genes} and is reproduced with permission from the author.}
    \label{fig:snp}
\end{figure}

Single-nucleotide polymorphisms (SNPs) are the most common type of genetic variation among human populations \cite{Shastry2009}. 
Analyzing SNPs holds significant value in many applied fields, as they can, among many other things, indicate a predisposition to various types of diseases, predict individual response to medication, and help trace human migration patterns (see Figure \ref{fig:snp}). 
A useful method for analyzing SNP data is through performing Principal Component Analysis (PCA), which reveals genetic variation or similarity among individuals and groups \cite{Abraham2014,Chen2016}.

\subsubsection*{High-Level Mathematical Description}
Suppose that the SNP data is given as a matrix $\mat{A} \in \R^{m \times n}$, where $m$ represents the number of individuals; $n$ represents specific SNPs across the genome; and
entry $\mat{A}_{ij}$ is the data for individual $i$ at SNP $j$, which can take values $0$, $1$ or $2$.
To apply PCA to SNP data, $\mat{A}$ must first be normalized:
\begin{equation}
\mat{B}_{ij}=\dfrac{\mat{A}_{ij} - \mu_{j}}{\sqrt{\dfrac{1}{2}\vct{\mu}_{j}\left(1 - \dfrac{1}{2}\vct{\mu}_{j}\right)}} \textrm{ , where } \vct{\mu}_{j} = \dfrac{1}{m}\sum_{i=1}^{m}\mat{A}_{ij}.
\end{equation}
After the normalized matrix $\mat{B}$ is obtained, the principal components can be extracted from the \textit{right singular vectors} of $\mat{B}$, which are obtained from computing the Singular Value Decomposition (SVD) of $\mat{B}$ \cite{Higham:blog:big6}.
The SVD of $\mat{B}$ is given by orthogonal matrices $\mat{U} \in \R^{m \times m}$, $\mat{V} \in \R^{n \times n}$, and diagonal matrix $\mat{S} \in \R^{m \times n}$, such that 
\begin{equation}
\mat{B} = \mat{U}\mat{S}\mat{V}^{^\intercal}.
\end{equation}
Given the SVD, the leading columns of $\mat{V}$ represent the principal components that capture the most variation in the data.
As such, a \textit{truncated} SVD, which aims to compute the leading $k \ll \min\{m, n\}$ singular triplets, is often more numerically efficient in this context.
The reduced representation, comprised of the leading $k$ principal components, can then be used to identify clusters or patterns in the dataset. 
In this context, performing a truncated SVD on $\mat{B} \in \R^{m \times n}$ can present a significant computational challenge, as $m$ and $n$ can be on the order of millions, depending on the particular input dataset \cite{ukbiobank,populussnp}.

\subsubsection*{The Potential Role of RNLA}

There are many practical approaches for finding a truncated Singular Value Decomposition \cite[e.g.,][]{EigenSVD,spectra_partial_svd,cusolver2025,magma_svd}. 
Among the methods that leverage modern computational strategies, Randomized Singular Value Decomposition (RandomizedSVD) and Randomized Subspace Iteration (RSI) stand out as some of the most widely used techniques in RNLA \cite{doi:10.1137/S0097539704442684,JMLR:v6:drineas05a,halko2011finding,RandLAPACK_rsvd}.

Despite their popularity, both RandomizedSVD and RSI are known to struggle with input matrices that have little variation in singular values (i.e., a slowly decaying spectrum).
This particular structure can be seen in many practical datasets, such as SNP data containing individuals from genetically distinct subpopulations, leading to highly structured genetic variation.
In such cases, an alternative yet less common RNLA method---Randomized Block Krylov Iteration (RBKI)---can be employed.
In comparisons of RBKI and RSI \cite{musco2015randomized}, RBKI can be substantially more accurate in certain problems, which has been substantiated by recent theoretical work and more realistic computational comparisons on an actual dataset
\cite{tropp2023randomized}.\footnote{The data set is sourced from \cite{flashpca}.}

\subsubsection*{Outlook and Open Challenges}

While  a solid theoretical foundation and compelling practical illustrations exist for justifying the use of RBKI \cite{tropp2023randomized}, a version of this method has yet to be implemented as part of high-performance linear algebraic software. 
A crucial aspect of implementing a practical RBKI-type method is ensuring its adaptability to large, sparse, and structured (Hermitian, positive semidefinite) inputs, enabling its application for (and beyond) SNP analysis in real-world settings.
Another important consideration is ensuring its accessibility to non-expert users by providing a simple interface for tuning algorithm parameters, as well as a quantitative measure of the error in each of the estimated singular triplets.
Importantly, the dominant computational cost of RBKI is attributed to level 3 BLAS operations (which involve interactions with the large input matrix); luckily, level 3 BLAS are the most efficient linear algebraic operations on modern systems. 
Thus, a variant of the RBKI-style algorithm is a promising candidate for inclusion as a driver-level method in the open-source RandLAPACK library \cite{randlapack}.

\subsection{Operator Inference} \label{application:oi}

Differential Equations (DEs) of the form
\begin{equation} \label{eqn:oi-de}
\frac{d\vct{x}}{dt} =  f(\vct{x}, \vct{\mu}),
\end{equation}
where $\vct{x}(t,\vct{\mu})$ is a state vector dependent on time, $t$, and system parameters, $\vct{\mu}$, are pervasive in scientific computing.  
For many practical situations, solving these DEs requires expensive numerical simulations. The expense of these simulations can be especially problematic in areas like control, optimization, and uncertainty quantification where such simulations must be run many times with different parameterizations to achieve a desired goal. One approach to reducing the expense of these numerical simulations is to approximate them with surrogate models, which can be run at a much lower computational cost \cite{kramer2024learning, farcas2023parametric,swischuk2020learningb}. 

\begin{figure}[!ht]
    \centering
    \includegraphics[width=0.7\linewidth]{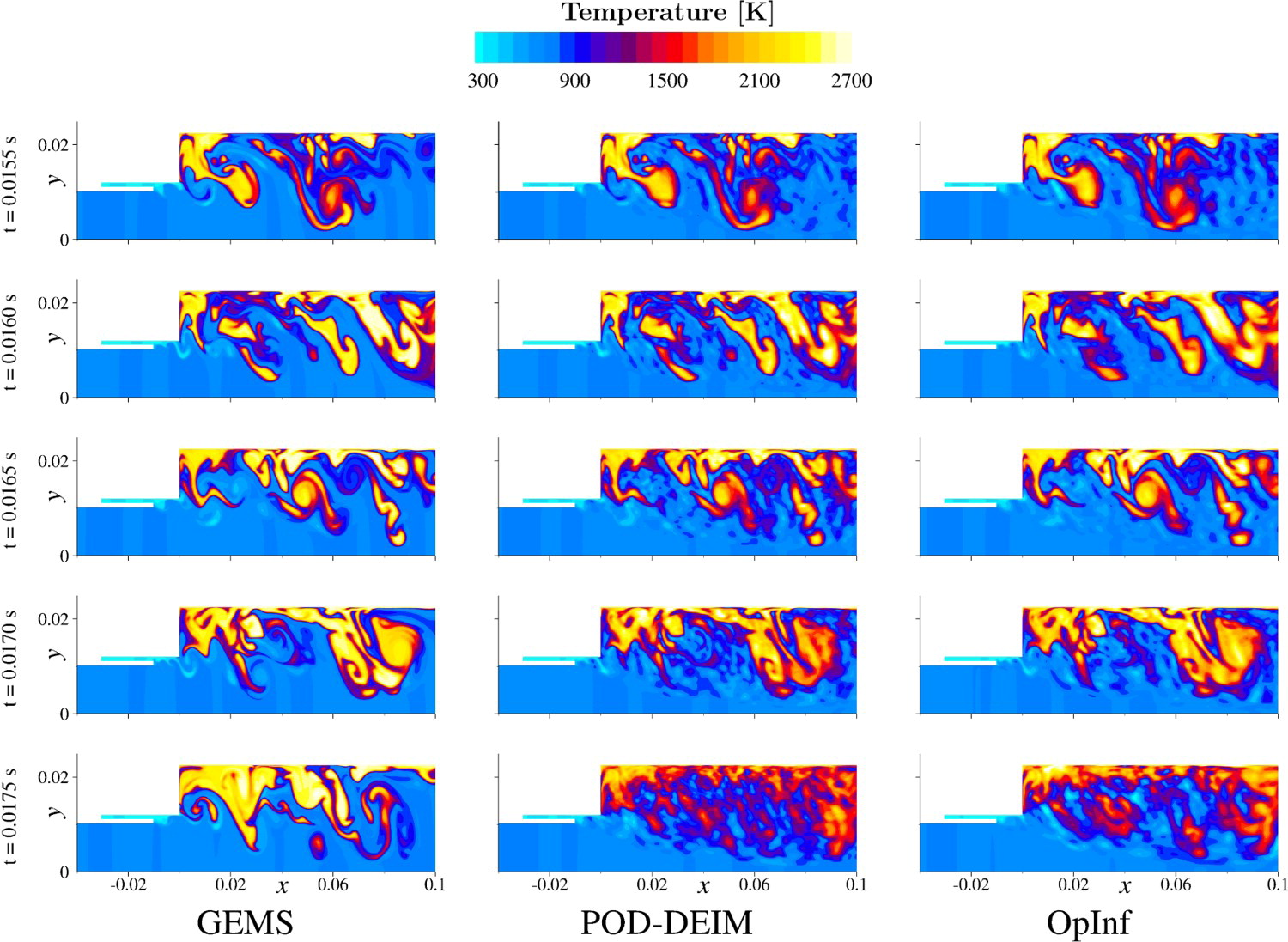}
    \caption{An example of the temperature field over time (paneled vertically) for a single-injector combustion problem. The GEMS column is generated from an expensive, high-fidelity simulation tool, while the Operator Inference (OpInf) column is generated from a low-fidelity model and captures the temperature field well, especially as time increases. This is a reproduction of Figure 7 from \cite{mcquarrie2021data} with permission from the authors.}
    \label{fig:oi}
\end{figure}

Operator Inference is one such surrogate modeling technique (see \cite{peherstorfer2016opinf} and Figure \ref{fig:oi}). Operator Inference fits a polynomial model to a set of reduced simulation results generated across numerical solutions of \eqref{eqn:oi-de} at different system parameters (and, optionally, different initial conditions) \cite{kramer2024learning}. 

\subsubsection*{High-Level Mathematical Description}
Operator Inference model has several steps.
First, the DE is numerically solved for $n$ distinct parameterizations, and the solutions at $m$ spatial points and $\tau$ time points are recorded for parametrization. The DE's numerical solutions, $\lbrace \vct{x}(t_{\ell}, \mu_j) : \ell=1,\ldots,\tau; j=1,\ldots,n \rbrace \subset \R^m$, are stored in a matrix, $\mat{X} \in \mathbb{R}^{m \times n\tau}$. Importantly, these DE solutions may already exist from previous numerical simulations and so we can directly make use of this data for Operator Inference.

Second, the matrix $\mat{X}$ is embedded into a lower dimensional space of dimension $k$ ($k < m$), resulting in a matrix $\hat{\mat{X}}$ \cite{kramer2024learning, farcas2023parametric,swischuk2020learningb}. 
The embedding of $\mat{X}$ is done using the leading $k$ left singular vectors of $\mat{X}$, denoted $\mat{U}_k \in \R^{m \times k}$, such that $\hat{\mat{X}} = \mat U_k^\intercal \mat{X}$.
Analogously, the embedded trajectories can be represented as $\lbrace \hat{\vct{x}}(t_\ell, \mu_j): \ell=1,\ldots,\tau;j=1,\ldots, n \rbrace \subset \R^k$.

The final stage of Operator Inference requires fitting a polynomial model---of order $p$---to the embedded numerical solutions. 
To represent the polynomial system, let $\vct{x}^h \in \R^{d_h}$ denote $h \in \mathbb{N}$ Kronecker products of $\vct{x} \in \R^m$ with itself and removing duplicate entries, which results in $d_h = \genfrac{(}{)}{0pt}{}{m+h-1}{h}$  \cite{kramer2024learning}.\footnote{As an example, if $\vct{x} \in \R^2$ has components $x_1$ and $x_2$, then $\vct{x}^2 = \vct{x} \otimes \vct{x}$ is three-dimensional with entries $x_1^2$, $x_1x_2$, and $x_2^2$.} 
Furthermore, for $h \in \mathbb{N}$, let $\mat A_h(\mu_j) \in \R^{k \times d_h}$ for $h=1,\ldots,p$ and $j=1,\ldots,n$.
Then, Operator Inference identifies $\mat A_h(\mu_j)$ for each $j$ by solving 
\begin{equation}
	\min_{\mat{A}_1(\vct{\mu}_j), \dots, \mat{A}_p(\vct{\mu}_j)} \sum_{\ell = 1}^\tau\left\| \sum_{h = 1}^p \mat{A}_h (\vct{\mu}_j) \hat{\vct{x}}^h(t_\ell, \mu_j) - \hat{\vct{x}}'(t_\ell, \mu_j)\right \|_F^2 + \sum_{h = 1}^p \lambda_h \|\mat A_h(\mu_j)\|_F^2,
\end{equation}
where $\hat{\vct{x}}'(t_\ell, \mu_j) \in \R^k$ represents a (possibly approximate) time derivative of the system at a time point $t_\ell$ and parameter $\mu_j$; and $\lbrace \lambda_h: h=1,\ldots, p\rbrace$  are non-negative regularization coefficients to prevent overfitting \cite{kramer2024learning}.

Once the polynomial model is determined, a dynamical system in $\R^k$ can be readily generated, and then multiplied by $\mat U_k$ to approximate the solution of the original DE in $\R^m$.

\subsubsection*{The Potential Role of RNLA}
In standard operator inference, computing $\hat{\mat{X}}$ from $\mat{X}$ can be expensive when $m$, $n$ and $\tau$ are moderately sized. RNLA can address this by efficiently computing low-rank approximations of the matrix $\mat{X}$. 
In fact, the RandomizedSVD \cite{halko2011finding} is already being used for this purpose \cite{kramer2024learning,benner2020operator, swischuk2020learningb,farcas2023parametric,kim2025physically}.  
RandomizedSVD has been shown to maintain similar accuracy to the full SVD \cite{rajaram2020randomized} and can speed-up the computation of the low-rank approximation thirty to fifty fold \cite{bach2019randomized}. 

\subsubsection*{Outlook and Open Challenges}
Recently, interpolative decompositions---that is, decompositions that select exact row subsets of a matrix---were shown to achieve an approximation error similar to that the RandomizedSVD \cite{osinsky2023close}. 
Randomized versions of these decompositions can potentially achieve even greater speed-ups compared to RandomizedSVD \cite{dong2023simpler}. 
To promote interpolative decompositions, theoretical developments must be made showing that such interpolative decompositions offer similar quality subspaces in terms of operator inference performance as the RandomizedSVD. Additionally, to improve the practicality of interpolative decompositions, such rank-adaptive approaches as those that exist for RandomizedSVD should be developed \cite{yu2018efficient,dong2023robust,pearce2025adaptive}; these approaches ensure the approximation satisfies a particular quality requirement.  In terms of software, providing a modular framework that allows for the easy substitution of different low-rank approximation techniques into Operator Inference implementations would go a long way to help with the exploration of the potential of such randomized techniques.

\subsection{Data Assimilation} \label{application:da}
Given a dynamical system, data assimilation (DA) aims to blend a computational model of a dynamical system with noisy observations of the system to estimate the system's states, parameters, or controls \cite{evensen2022data,asch2016book}.
%
%
The Strong Constraint 4D-Variational approach (SC-4DVAR) is a well-known DA framework used to estimate the initial condition of a system under the assumption that the computational model of the dynamical system is accurate \cite{Kalnay_B2003}.

\subsubsection*{High-Level Mathematical Description}
Data assimilation for state estimation combines information from three different sources---the prior best estimate of the initial state of the system, the model of the system, and observations of the system---to produce an improved estimate of the initial state of the dynamical system.

Formally, let us denote the unknown state at time $t_0$ as $\vct{x}^{\rm t}_0 \in \R^n$ and the prior best estimate by 
$\vct{x}^b_0 \in \R^n$.
The computational model---here, a set of differential equations---evolves an initial state $\vct{x}_0 \in \mathbb{R}^n$ at time $t_0$ to future state $\vct{x}_i \in \mathbb{R}^n$ according to
\begin{equation}\label{eqn:model}
    \vct{x}_i = \mat{M}_{t_0 \rightarrow t_i} (\vct{x}_0)\,, \text{ where } \mat{M}_{t_0 \rightarrow t_i} : \R^n \to \R^n, \> 1 \leq i \leq \tau.
\end{equation}

The observations, $\vct{y}_i \in \R^{n_{\rm obs}}$, of the dynamical system are taken at times $t_i$, $1 \leq i \leq \tau$, and are related to the computational model by 
\begin{align}
    \vct{y}_i = \mat{O}(\vct{x}^{\rm t}_i) + \varepsilon_i, 
\end{align}
where $\mat O:\R^n \to \R^{n_{\rm obs}}$ is the observation operator;%
\footnote{Typically, $n_{\rm obs} \ll n$.}
and $\varepsilon_i$ represents observation error.

SC-4DVAR computes an estimate of $\vct{x}_{0}^{\rm t}$ as a solution to
\begin{equation}\label{eqn:4dvar}
    \begin{aligned}
        \min_{\vct{x}_0\in \R^n} & \frac{1}{2}  \left(\vct{x}_0 - \vct{x}^b_0\right)^\intercal \vct\Gamma_p^{-1} \left(\vct{x}_0 - \vct{x}^b_0\right) 
     + \frac{1}{2} \sum_{i=1}^{\tau} \left( \mat{O}(\vct{x}_i) - \vct{y}_i\right)^\intercal \mat R_i^{-1} \left( \mat{O}(\vct{x}_i) - \vct{y}_i\right), \\
        \text{subject to: } & \vct{x}_i = \mat{M}_{t_0 \rightarrow t_i} (\vct{x}_0)\,,\quad  1 \leq i \leq \tau,
    \end{aligned}
\end{equation}
where $\mat{\Gamma}_p$ and $\mat{R}_i$ are the variances of the initial estimate, $\vct{x}_0^{\rm b}$ and $\vct{\varepsilon}_i$, respectively.

\subsubsection*{The Potential Role of RNLA}
Minimizing the cost function in \eqref{eqn:4dvar} using quasi-Newton methods or Newton's method requires solving multiple linear systems involving the Hessian (or its approximation) of the cost function. 
As the typical size of $\vct{x}_0$ is in billions and that of $\vct{y}$ is in tens to hundreds of millions, it is impractical to form the Hessian explicitly, and the resulting linear system must be solved using iterative methods that use the Hessian implicitly. 
Furthermore, the linear systems involving the Hessian can often be ill-conditioned, which necessitates efficient preconditioners. 

RNLA offers a relatively inexpensive approach to generating such preconditioners through RandomizedSVD, Nystrom approximations and single view approaches \cite{halko2011finding, tropp2017practical, martinsson2020randomized}. 
By using these preconditioners, the minimization of \eqref{eqn:4dvar} can be meaningfully accelerated, and the resulting procedure enjoys theoretical guarantees \cite{subrahmanya2024randomized}.

\subsubsection*{Outlook and Open Challenges}
The primary challenges are to create more sophisticated randomized techniques that account for other DA methods beyond SC-4DVAR; to address Hessians or their approximations that have slowly decaying spectra; and incorporating randomized algorithms for problems that use multi-fidelity models.

\subsection{Lattice Quantum Chromodynamics I: Operator Simulation} \label{application:qcd_overlap}
Quantum chromodynamics (QCD) is the area of theoretical particle physics that studies the strong interaction between quarks and gluons. 
At low energies (i.e., for hadrons), lattice QCD is a formulation of this theory in which space-time is discretized and simulated on a four-dimensional lattice. 
Lattice QCD specifies dynamics that are too complicated to analyze analytically,
yet these dynamics can be simulated numerically allowing us to predict important 
physical quantities (e.g., hadron masses, form factors, etc.). 

\subsubsection*{High-Level Mathematical Description}
The strong interaction is described by the Dirac equation, $(\mathcal{D} + m) \vct{\psi}(\vct{x}) = \vct{\eta}(\vct{x})$, where  
$\vct x$ is a position in time and three-dimensional space;
$m$ is a scalar parameter; 
$\vct{\psi}(\vct x), \vct{\eta}(\vct x) \in \C^{12}$ represent quark fields;%
\footnote{The 12 variables correspond to all possible combinations of three colors and four spins.}
and $\mathcal{D}$ is the {Dirac operator} 
\begin{equation}
    \mathcal{D} = \sum_{i = 0}^3 \mat{\Gamma}_i \otimes \left(\frac{\partial}{\partial {\vct x}_i} + \mat{A}_i(\vct{x})\right),
\end{equation}
where $\mat{\Gamma}_i \in \C^{4 \times 4}$ is a generator of the Clifford algebra; $\partial/\partial \vct{x}_i$ is a derivative operator on the ``color space'' of the quarks; and $\mat{A}_i({\vct x}) \in \C^{3 \times 3}$ is an element of the Lie algebra $\mathfrak{su}(3)$  for $i=0,1,2,3$  \cite{dira58}. 

For simulations in lattice QCD, the Dirac equation is discretized on a toroidal lattice with $N_t$ lattice points in time and $N_s$ lattice points in each of the three spatial dimensions.
For example, the popular Wilson-Dirac discretization \cite{wils77} results in an operator, $\mathcal{D}_W$, where 
$(\mathcal{D}_W \vct \psi)(\vct x)$ is 
\begin{equation} \label{eqn:wilson}
\begin{aligned}
\frac{m_0 + 4}{a}\vct\psi(\vct{x}) - \frac{1}{2a}\sum_{i = 0}^3 \left[ \left((\mat I_4 - \mat\Gamma_i) \otimes \mat U_i(\vct x)\right)\vct\psi(\vct x + a\vct{e}_i) + \left((\mat I_4 + \mat \Gamma_i) \otimes \mat U_i\herm(\vct x-a\vct{e}_i)\right)\vct\psi(\vct x - a\vct{e}_i)\right],
\end{aligned}
\end{equation}
where 
$m_0$ determines the quark mass;
$a$ is the distance between lattice points in the discretization;
$\mat I_4$ is the $4 \times 4$ identity matrix;
the gauge links, $\mat U_i(x) \in \mathbb{C}^{3 \times 3}$, are elements of the Lie group $\text{SU}(3)$;
$\cdot\herm$ indicates the Hermitian transpose of a matrix;
and $\vct e_i$ are the standard basis elements for $i=0,1,2,3$.
The Wilson-Dirac discretization can be simulated by solving a sequence of systems of equations with thousands 
of linear systems, which can be done efficiently by multigrid and domain decomposition methods \cite{frommer2014adaptive,luscher2004solution}. 

\begin{figure}
    \centering
    \includegraphics[width=0.9\linewidth]{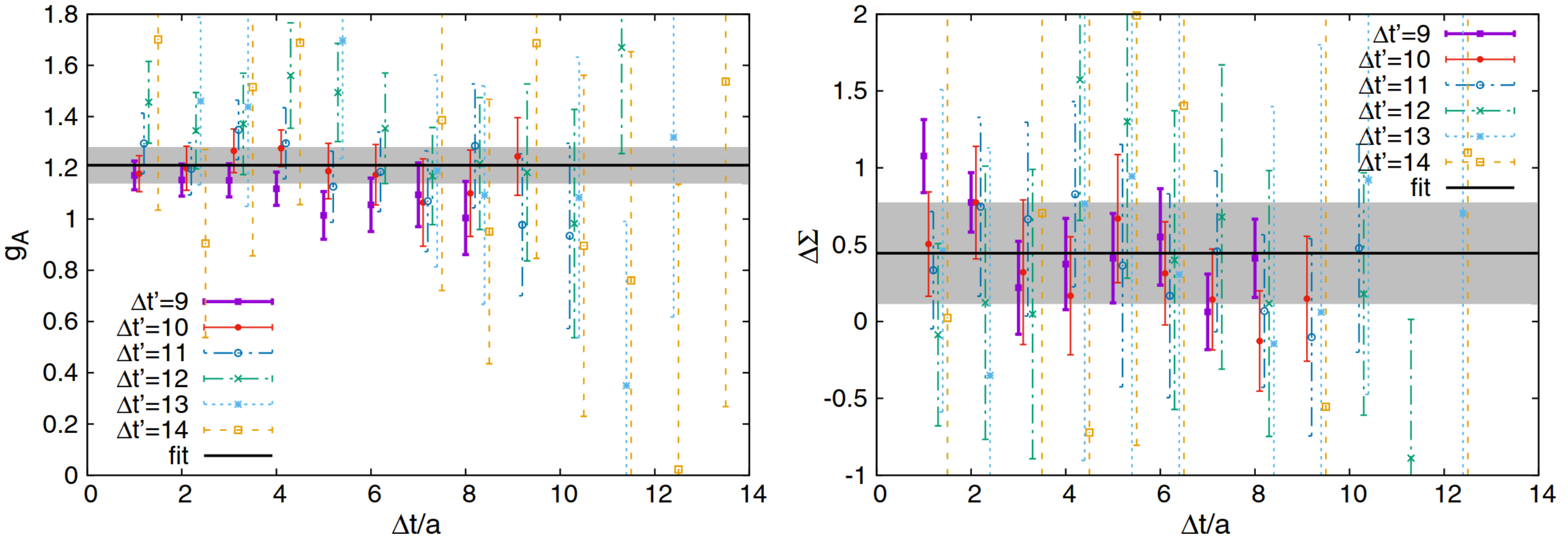}
    \caption{Analysis of nucleon isovector axial charge $g_A$ (left) and nucleon singlet axial charge $\Delta\Sigma$ (right) using \emph{overlap fermions} at lattice spacing $a$ and different time separations $\Delta t$. Because of the exact preservation of chiral symmetry on the lattice, charges can be extracted without spurious operator mixing. The horizontal line and the grey band show the fitted value and statistical error of each charge, yielding $g_A = 1.210 \pm 0.072$ and $\Delta\Sigma = 0.44 \pm 0.33$. The horizontal plateaus indicate that excited‑state contaminations have decayed away and that the ratio defining the effective charge has reached a true ground‑state. This is Figure 9 from \cite{yamanaka2018nucleon} and is reproduced under CC-BY 4.0.}
    \label{fig:qcd}
\end{figure}

An alternative discretization of the Dirac equation that preserves chiral symmetry (unlike $\mathcal{D}_W$),
which is essential for the simulation of such observables as nucleon axial charges (c.f., Figure \ref{fig:qcd}), is the Neuberger overlap operator \cite{neub98}
\begin{equation}\label{eq:overlap_operator}
\mathcal{D}_N = \rho \mat I_p + \mat{\Gamma}_5\text{sign}(\mat{\Gamma}_5 \mathcal{D}_W),
\end{equation}
where $\rho > 1$ is a regularization term;
$p = 12N_t N_s^3$;
$\mat I_p$ is the $p \times p$ identity matrix;
$\mat \Gamma_5 = \mat I_{N_t N_s^3} \otimes (\mat \Gamma_0 \mat \Gamma_1 \mat \Gamma_2 \mat \Gamma_3 ) \otimes \mat I_3$;
and $\text{sign}$ is the matrix sign function \cite[Ch. 5]{high08}. 
In the presence of a {nonzero chemical potential}, the matrix resulting from evaluating $\mat \Gamma_5 \mathcal{D}_W$ is non-Hermitian.
The Neuberger overlap operator can be simulated by solving a sequence of systems of linear equations, typically done by using iterative methods that only require matrix-vector products (e.g., Krylov subspace methods). 

\subsubsection*{The Potential Role of RNLA}

Solving the systems of linear equations for simulating with the Neuberger overlap operator is challenging for multiple reasons. 
First, each system of linear equations is large: the minimum number of lattice points in time and space for realistic simulations are $N_t = 64$ and $N_s = 32$, which results in square systems of dimension  $12 \times 64 \times 32^3 = 25, 165, 824$ by $25, 165, 824$. 
Second, each system requires evaluating $\mathcal{D}_W$, which can only be accessed through matrix-vector products (i.e., given $\vct{v}$, we can evaluate $\mathcal{D}_W \vct v$) in product-level code. 
Consequently, approximating the system's $\text{sign}(\mat \Gamma_5 \mathcal D_W)$ term is done via
simple polynomial Krylov subspace methods built upon matrix-vector products, which converge very slowly for the sign function. 
Finally, if the system's $\mat \Gamma_5 \mathcal D_W$ is non-Hermitian, then the orthogonalization cost in such Krylov methods grows quadratically in the number of iterations, often limiting the maximum number of iterations (and thus accuracy) that can be reached. 

RNLA methods have been recently introduced that dramatically reduce the orthogonalization cost. 
Such RNLA methods are of two flavors: 
a randomized Gram--Schmidt process that performs orthogonalization using much smaller, sketched basis vectors \cite{balabanov2022randomized,cortinovis2024speeding};
or a truncated orthogonalization followed by a ``basis whitening'' step, which retrospectively makes the sketched basis  orthogonal \cite{GuettelSchweitzer2023,palitta2025sketched}. 
With either RNLA method, the cost of orthogonalization scales linearly in the number of iterations, significantly reducing the cost of approximating $\text{sign}(\mat \Gamma_5 \mathcal{D}_W)$. 

Furthermore, RNLA methods are well-suited for recycling information to speed up computations of a sequence of matrix function problems \cite{burke2024krylov}. Such a sequence of matrix function problems arise when simulating with the Neuberger overlap operator, which would further accrue the benefits of RNLA methods in lattice QCD simulations. 

\subsubsection*{Outlook and Open Challenges}

While randomized Krylov methods show very promising experimental performance in exploratory studies (also for small-scale lattice QCD problems, see \cite[Section 5.3]{GuettelSchweitzer2023}), their theory is insufficient for lattice QCD problems. While their behavior is at least partially understood for entire functions, much less is known for functions with singularities (such as the matrix sign function), and it cannot be ruled out yet that the methods can completely break down. Another problem preventing the adoption of randomized Krylov methods in production level code is the lack of reliable and cheaply computable stopping criteria, a topic that has not received much attention yet.%
\footnote{Some progress has been made for randomized row and column projection methods \cite{pritchard_towards_2023,pritchard_solving_2024}, but not randomized Kyrlov methods.}

\subsection{Lattice Quantum Chromodynamics II: Implicit Trace Estimation} \label{application:qcd_trace}
Lattice QCD discretizations can be used to compute physical observables such as disconnected quark loops \cite{durr2012sigma}, the Dirac spectrum \cite{cossu2016stochastic}, and fluctuations of topological charge \cite{luscher1982topology}. 

\subsubsection*{High-Level Mathematical Description}
For example, $\tr(\mathcal D_W^{-1})$, where $\mathcal D_W$ is the Wilson discretization of the Dirac operator given in \eqref{eqn:wilson}, indicates quark-loop contributions.
However, computing $\tr(\mathcal D_W^{-1})$ is challenging owing to the size of $\mathcal D_W$, which 
has dimension $12N_t N_s^3 \times 12 N_t N_s^3$, where $N_t$ is the number of discretization points in time and $N_s$ is the number of discretization points in each spatial dimension. 

\subsubsection*{The Potential Role of RNLA}

RNLA methods can approximate $\tr(\mathcal D_W^{-1})$ via
\begin{equation}
\tr\left( \mathcal{D}_W^{-1} \right)
= 
\mathbb{E}\left[\vct{z}^\ast \mathcal{D}_W^{-1}\vct{z}\right]
\approx 
\frac1s\sum_{j=1}^s \vct{z}_j^\ast \mathcal{D}_W^{-1}\,\vct{z}_j,
\end{equation}
where each $\vct{z}_j\in\C^{12N_tN_s^3}$ is a random ``probe'' vector with $\mathbb{E}[\vct{z}_j \vct{z}_j\herm] = \mat{I}_{12N_tN_s^3}$. 
These RNLA trace estimators replace computing the inverse with solving $s$---the number of samples---linear systems (i.e., $\mathcal{D}_W^{-1} z_j$), which can be done efficiently by  multigrid solvers or multigrid-preconditioned Krylov methods \cite{frommer2014adaptive,luscher2004solution}.
Thus, such randomized trace estimators are efficient when the sample size, $s$, is much smaller than $12N_tN_s^3$. 

In one avenue, randomized trace estimation advancements aim to reduce the number of samples using several strategies such as (inexactly) deflating low eigenmodes \cite{romero2020multigrid}; or constructing lattice-structure-aware probing vectors \cite{stathopoulos2013hierarchical,foley2005practical}. 
In another avenue, randomized trace estimation advancements use {multi-level estimators}. Multi-level estimators introduce a hierarchy of samples evaluated at different accuracies to balance variance reduction of the estimator and computational cost. A attractive multi-level estimator couples the sample hierarchy to the multigrid hierarchy that is used for solving the linear systems.

\subsubsection*{Outlook and Open Challenges}
One particular open problem is how to best quantify the trace estimation error in order to decide how many samples to use or when to stop iterating. Current practice mostly relies on heuristics (e.g., the empirical variance of the $s$ samples). Rigorous, \emph{a‑posteriori} error bounds---perhaps via concentration inequalities tailored to the lattice structure---would allow certificates of accuracy at fixed cost.
Another possible area of improvement is the adaptive selection of probing vectors. While hierarchical probing uses a fixed coloring of the lattice, adaptive strategies that ``learn'' which patterns of the gauge field contribute most to variance---and concentrate probes there---could further reduce computational cost.
Finally, embedding these RNLA trace estimators into production lattice frameworks like Chroma \cite{edwards2005chroma}, Grid \cite{boyle2016grid} or QUDA \cite{clark2010solving} with user‑transparent interfaces and auto‑tuning of $s$, precision, and probe structure remains an essential step toward wider adoption.

\pagebreak 

\section{Conclusion} \label{section:conclusion}
The applications presented herein show that RNLA is maturing from a collection of theoretical techniques into a fundamental computational paradigm for modern scientific computing.  Across the discussed domains, RNLA methods have proven capable of addressing computational challenges that classical numerical linear algebra cannot handle at the scales demanded by current and future scientific applications.
To ensure its continued utility and wide-spread adoption across scientific domains, RNLA must address several critical areas.

\subsubsection*{Algorithmic foundations for structure-aware methods} The development of principled frameworks for exploiting problem-specific structure represents a major challenge. Future research must develop approaches that can systematically identify and exploit mathematical structure through randomization. 

\begin{enumerate}
    \item Adaptive structure detection: methods are needed to automatically identify exploitable structure and adapt randomization strategies accordingly (e.g., correlation structure as identified in \S \ref{application:gwas2}).
    \item Structure-preserving guarantees: methods require error guarantees that ensure preservation of domain-specific properties such as non-negativity or physical constraints
    (e.g., \S \ref{application:hyperspectralmixing}).
\end{enumerate}

\subsubsection*{Hardware-Algorithm Co-Design} The computational advantages of RNLA can only be fully realized when exploiting modern computer architectures. The applications discussed herein consistently identify GPU acceleration and mixed-precision as a promising avenue.

\begin{enumerate}
    \setcounter{enumi}{2}
    \item Memory hierarchy optimization: data movement costs can constitute a significant portion of computational time, particularly when handling large datasets (e.g., \S\ref{application:ct}). Thus, there is a need to develop cache-aware randomized algorithms that minimize data transfer overhead. This optimization becomes especially critical in GPU-accelerated computing, where data transfer between host and device memory can dominate overall execution time.
    \item Precision-aware randomization: mixed-precision strategies have the potential to vastly decrease memory requirements. However, the way that reduced or mixed precision affect the performance of RNLA algorithms is poorly understood. Research is needed to develop randomized algorithms that can adaptively adjust precision while maintaining theoretical guarantees.
\end{enumerate}

\subsubsection*{Software Infrastructure for Reproducibility}
The transition from proof-of-concept to production-ready software represents a critical bottleneck in RNLA adoption. While established libraries exist for linear algebra, equivalent specifications and implementations of RNLA are emerging \cite{murray_randomized_2023, patel_rlinearalgebra_2025}.

\begin{enumerate}
    \setcounter{enumi}{4}
    \item Modular frameworks: the diverse requirements across applications demand modular software frameworks that allow domain scientists to easily substitute different randomized components. The RandLAPACK initiative mentioned in \S \ref{application:snp} represents an important step in this direction. Other initiatives such as RLinearAlgebra emphasize the role of composable algorithms \cite{patel_rlinearalgebra_2025}.
    \item Reproducibility: ensuring consistent results across multiple runs presents fundamental challenges for RNLA deployment. This will require developing systematic approaches to variance control mechanisms that bound resulting variability within acceptable tolerances, and creating standardized protocols for comparing the performance of different randomized algorithms.
    \item Integration with scientific software ecosystem: each application domain has established data formats, software ecosystems, and computational workflows. RNLA libraries must integrate seamlessly into these existing frameworks rather than requiring adoption of new software stacks.
\end{enumerate}

\subsubsection*{Theoretical foundations for practical development} The gap between theoretical guarantees and practical performance represents a fundamental challenge for RNLA. Applications require more computationally tractable error estimation and stronger theoretical guarantees.

\begin{enumerate}
    \setcounter{enumi}{7}
    \item Practical error control: while most methods provide probabilistic guarantees, the gap between error bounds and practical error control remains significant. Developing theory that bridges probabilistic guarantees with actionable error estimation is essential for the adoption of RNLA methods by the scientific community.\footnote{Some works that focus on tractable error estimation include \cite{tropp2017practical,pritchard_towards_2023,pritchard_solving_2024}.}
\end{enumerate}

By addressing these eight, application-motivated challenges, RNLA will become an essential component of a plethora of large-scale scientific computational workflows.

\newpage

\printbibliography

\end{document}